\documentclass[12pt]{article}
\usepackage{dbw,amsmath,amssymb,amsfonts,psfig,picinpar,xspace}
\usepackage[round]{natbib}

\newcommand{\sref}[1]{\S~\ref{sec:#1}}

\newcommand{\fref}[1]{Figure~\ref{fig:#1}}

\newenvironment{code}{\begin{list}{\hspace*{0em}}{}}{\end{list}}
\newcommand{\old}[1]{}
\newcommand{\rmap}{\ensuremath{\text{RandomMap()}}}
\newcommand{\While}{{\tt while\ }}
\newcommand{\Repeat}{{\tt repeat\ }}
\newcommand{\Until}{{\tt until\ }}
\newcommand{\For}{{\tt for\ }}
\newcommand{\To}{{\tt to\:\,}}
\newcommand{\DownTo}{{\tt downto\:\,}}
\newcommand{\If}{{\tt if\ }}
\newcommand{\Not}{{\tt not\ }}

\newcommand{\Return}{{\tt return\ }}
\newcommand{\Output}{{\tt output\ }}

\newcommand{\s}{\hspace*{3ex}}
\newcommand{\sspace}{$\langle\text{state space}\rangle$}

\newcommand{\R}{{\mathbb R}}

\begin{document}
\title{\vspace*{-50pt}
       How to Couple from the Past\\ Using a Read-Once Source of Randomness}
\author{\vspace*{-3pt}
David Bruce Wilson\\Microsoft Research}
\date{}
\maketitle

\vspace*{-20pt}
\begin{abstract}
We give a new method for generating perfectly random samples from the
stationary distribution of a Markov chain.  The method is related to
coupling from the past (CFTP), but only runs the Markov chain forwards
in time, and never restarts it at previous times in the past.  The
method is also related to an idea known as PASTA (Poisson arrivals see
time averages) in the operations research literature.  Because the new
algorithm can be run using a read-once stream of randomness, we call
it read-once CFTP.  The memory and time requirements of read-once CFTP
are on par with the requirements of the usual form of CFTP, and for a
variety of applications the requirements may be noticeably less.  Some
perfect sampling algorithms for point processes are based on an
extension of CFTP known as coupling into and from the past; for
completeness, we give a read-once version of coupling into and from
the past, but it remains unpractical.  For these point process
applications, we give an alternative coupling method with which
read-once CFTP may be efficiently used.
\end{abstract}

\section{Introduction}

One of the mantras of ``coupling from the past'' (CFTP), a class of
algorithms for generating perfectly random samples from a Markov
chain, is that one needs to be prepared to re-use old random coins.
This would appear to rule out any possibility of running CFTP with a
read-once stream of random coins, such as a Geiger counter,
thermal noise (used by the Intel hardware random number generator
\citep{jun-kocher:intel-rng}), or other
physical source of truly random coins, short of storing all the random
values somewhere.  Nonetheless we give here a simple variation on
CFTP, whose time and memory usage is competitive with the current
prevalent version of CFTP, but which outputs perfectly random samples
using just a read-once source of randomness.  Even with re-readable
sources of pseudorandom coins, which come with seeds that allow one to
regenerate previously generated values, there can be advantages to
using the read-once version of CFTP, particularly if many independent
samples are desired.

We give a more detailed review of CFTP in a later section, but for now
we state that it is a method for generating random samples from the
steady state distribution of a Markov chain, when the Markov chain is
implemented by repeatedly applying a randomizing operation to a point
in the state space.  The method is based on the principle that a
Markov chain that has already been running for an infinitely long time
has already reached its stationary distribution.  To obtain a random
sample, CFTP ``figures out'' what state the Markov chain is in at a
given time, by looking at a finite but unbounded number of randomizing
operations used prior to that time.  Usually the ``figuring out'' part
requires cleverness on the part of the algorithm designer, and
different techniques are used in different applications.  Rather than
extend the range of applications to which CFTP may be applied, our
purpose here is to give a variation on the method which may be used
with most of these applications.

\medskip
\noindent{{\bf Main Result:\ }\it
Every version of CFTP for which
\begin{enumerate}
\item i.i.d.\ randomizing operations are used to do the updates,
\item the algorithm produces the random sample, in its entirety
      and with full precision, after composing finitely many random maps, and
\item the random maps can be evaluated at a given starting state
      without affecting coalescence detection,
\end{enumerate}
can be done with a read-once stream of random coins.
Furthermore, the expected running time and memory usage are never worse
by more than a (small) constant factor.
}

The three conditions of the main result are satisfied by most algorithms
that one would
normally think of as CFTP, with just a couple of exceptions.  The
exceptions to condition 1 are a few algorithms, which might be more
properly described as ``coupling {\it into and from\/} the past,''
which use a separate Markov chain running backwards in time, rather
than an i.i.d.\ process, to generate the randomizing operations used
to do the updates.  
The principal exceptions to condition 2 are
algorithms given by \cite{berg-steif:codings}, and
\cite{haggstrom-steif:finitary} for infinite spin systems, where there
is no hope of outputting a sample in finite time, but where there is a
``virtual infinite sample,'' any part of which can be revealed to
someone who asks to see it.  For these algorithms, a re-readable
source of randomness appears to still be required.
\cite{moller:conditional} gave an algorithm for the autogamma
distribution which outputs not the sample but a neighborhood
containing the sample after composing finitely many maps,
so it too does not satisfy condition 2.
But \cite{wilson:multishift} gave a modification which is not only
faster but also satisfies the three conditions of the main result,
thereby allowing us to use read-once CFTP.
At present there are no exceptions to condition 3.

One advantage of
read-once CFTP over the prevalent version of CFTP is that one does not
need to keep track of pseudorandom number generator seeds.  As CFTP is
typically currently used, for even one sample the program keeps track
of seeds for a number of independent streams of pseudorandom numbers.
When many independent samples are desired, many independent streams
are required.  This independence requirement could be a problem if one
is using e.g.\ the pseudorandom number generator that comes standard
with Unix (BSD or libc5), where even if one believes that the stream of
numbers produced from any given seed is adequately random, the streams
produced using different seeds are quite decidedly {\it not\/}
independent.  (The streams started by different seeds are correlated
to an extent that is quite shocking to someone expecting independence.)
Using read-once CFTP for even a large number of samples,
only one good-quality stream of pseudorandom numbers is needed.

Read-once CFTP is also advantageous in situations where
storage is currently used for each time step, even when there is a
re-readable source of random coins.  In some cases it is not feasible
to generate an entire random map at once, so the algorithm instead
maintains partial information about each random map, which is then
updated each time the random maps are revisited; examples are given by
\cite{lund-wilson:storage} and \cite*{mira-moller-roberts:slice}.  Read-once
CFTP never revisits a random map it has seen before, so it is not
necessary to either store this partial information, or, more
importantly, to write code to update this partial information.  Many
point process algorithms (see e.g.\ 
\cite{kendall-moller:exact-spatial}) also store information for each
time step.  Although these point process algorithms often do not
satisfy condition 1 of the main result, in \sref{point-process} we are
still able to use read-once CFTP to sample from these point processes,
thereby reducing the storage requirements.

The read-once version of CFTP given here will satisfy additional
pleasant run-time properties not mentioned in the claim.  In some
cases, such as applications of CFTP to Bayesian inference, the
read-once version of CFTP may be noticeably faster.  Other run-time
properties would be appreciated by someone concerned about the
sociological phenomenon of an impatient user introducing bias by
aborting and restarting the algorithm.  For instance, the distribution
of running times will have exponentially decaying tails, assuming that
the effort to apply a single random map does not itself have a run
time distribution with fat tails.  The usual version of CFTP also has
this property, but the read-once version of CFTP has another
favorable run time characteristic not shared by the usual version of
CFTP.

As \cite{fill:interruptible} has pointed out, the usual version of
CFTP will on occasion enter a state where the expected additional
running time before outputing an answer can be large.
When this happens,
the user may be tempted to abort and start over.
In contrast, under conventional assumptions (explained in
\sref{performance}) about the underlying random map procedure, the
read-once version of CFTP given below does not have this property.
For it the
expected time to completion is never larger than it would be if the user
aborted and started over.  A stochastic domination version of this
statement also holds, so it should be the case that the user is never
tempted to abort.  Thus we could say that the algorithm is
``temptation free.''
Despite the algorithm being temptation-free, the
user with a specific deadline (and re-readable randomness) may still
prefer to use an ``interruptible'' algorithm such as Fill's algorithm
\citep{fill:interruptible}.

In the remainder of this article we review CFTP in \sref{cftp}, and
then give two derivations of read-once CFTP, the first one (in
\sref{rocftp}) starts from CFTP, and the second (in \sref{pasta})
starts from another idea known as PASTA (Poisson arrivals see time
averages).  In \sref{performance} we characterize the performance of
read-once CFTP.  Many interesting applications of CFTP are to
unbounded state spaces, and in \sref{unbounded} we give a variation of
a subroutine of
read-once CFTP that makes it easier to use in these contexts.  In
\sref{ciaftp} we review the coupling into and from the past (CIAFTP)
algorithms, which do not satisfy the first condition (independence of
random maps) required by read-once CFTP.  We give a read-once version
of CIAFTP in \sref{rociaftp}, but it is not very satisfying.  As the
principal applications of CIAFTP are point processes, we explain in
\sref{point-process} how to sample from these point processes using
instead the version of read-once CFTP in \sref{unbounded}.

\section{Background on coupling from the past}
\label{sec:cftp}

Before describing the read-once version of CFTP, we first review the
usual version of CFTP.  More expanded explanations are given by
\cite{propp-wilson:exact-sampling}, \cite{fill:interruptible},
\cite{propp-wilson:cftp-aug}, and \cite{wilson:multishift}.

CFTP requires a randomizing operation which preserves the probability
distribution $\pi$ from which we wish to sample.  There are many maps
from the state space to itself; the randomizing operation effectively
picks a random such map according to some distribution.  Let us
consider a toy example: suppose $\pi$ is the uniform distribution on
the state space of permutations on $n$ letters.  One possible
randomizing operation would pick a random number $i$ between $1$ and
$n-1$, and then flip a coin $c$ to decide whether to rearrange the
items in positions $i$ and $i+1$ so that they are in sorted order or
in reverse-sorted order.  If we perform this operation on a uniformly
random permutation, the result will also be a uniformly random
permutation, so we say that the randomizing operation preserves the
uniform distribution $\pi$.  The (random) pair $(i,c)$ may be used to
update any given permutation, so it represents a (random) function or
map from the state space to itself.  We obtain a Markov chain by
applying the randomizing operation over and over again to a given
state; different randomizing operations may give rise to the same
Markov chain.

We assume that the randomizing operation is given to us as the
procedure RandomMap().  Each time that RandomMap() is called, it
returns some representation of a random map (such as a random pair
$(i,c)$ in the above example), and the random map is independent of
all random maps previously generated.  Let $f_{-t}$ denote the map
returned the $t^{\text{th}}$ time RandomMap() is called, which we view
as the randomizing operation that occured at time $-t$.  If a Markov
chain is in state $x$ at time $-t$, then at time $-t+1$ it will be in
state $f_{-t}(x)$.  Thus we view the randomizing operations as having
been started infinitely far in the past, and they run up until time
$0$.  Let $F_{-t}$ denote the composition of $f_{-1} \circ \cdots
\circ f_{-t}$, i.e.\ the net effect of the $t$ randomizing operations
prior to time $0$.  If we somehow obtained a random state $x$
distributed according to $\pi$, then since the randomizing operation
preserves $\pi$, $F_{-t}(x)$ will also be distributed according to
$\pi$.

It is easy to see that the event that there is some $t$ such that
$F_{-t}$ maps the state space to one value, occurs with probability
either $0$ or $1$.  Usually it is not hard to ensure that this
probability is positive, so let us assume that the probability is $1$.
If $F_{-t}$ maps the state space to a single value, then for any
$t'>t$, $F_{-t'}$ will also map the state space to this same value.
So with probability $1$, all but finitely many of the random variables
$F_{-1}(x), F_{-2}(x), F_{-3}(x), \dots$ will take the same.  Since
this common value is independent of $x$, for convenience we denote it
by $F_{-\infty}(*)$.  Since the random variables $F_{-1}(x),
F_{-2}(x), F_{-3}(x), \dots$ are each distributed according to $\pi$,
and with probability $1$ they converge to the random variable
$F_{-\infty}(*)$, this random variable must also be distributed
according to $\pi$.

CFTP, which is expressed abstractly as in \fref{citp}, works by
determining and then outputting the random variable $F_{-\infty}(*)$.
Either CFTP runs forever with probability $1$, or else with
probability $1$ it successfully determines the state $F_{-\infty}(*)$ of
the Markov chain at time $0$, which is distributed exactly according
to the desired distribution $\pi$.

\begin{figure}[phtb]
\begin{code}
\item  $F := \langle\text{identity map}\rangle$
\item  \While \Not Singleton(ImageOf($F$))
\item  \s $F := F \circ \rmap$
\item  \Return ElementContainedIn(ImageOf($F$))
\end{code}
\caption{
High level pseudocode for coupling from the past.  Since random maps
are composed going backwards in time, the algorithm might be more
properly called coupling {\it into\/} the past.
It has been observed many times that reversing the order of
composition in the third line would result in a biased algorithm.
}
\label{fig:citp}
\end{figure}

The fact that composing maps backwards in time gives information about
the state at time $0$, which is then a perfectly random sample, appears
to have been first noted and exploited by \cite{letac:contraction}.
\cite{diaconis-freedman:rfuncs} give a survey of this and related
work.  The main use for which this principle was used was to prove
the existence of stationary distributions of Markov chains.
Algorithms based on this principle for sampling
from nontrivial distributions weren't developed until many years
later.  The basic problem was a lack of effective means of determining
when to stop composing the maps.  The first (nontrivial) algorithms
based on the ``state at time zero is random'' principle was a random spanning
tree algorithm due to \cite{broder:tree} and \cite{aldous:tree}, and
the dead-leaves process (see \citep{jeulin:dead-leaves}).
The tree algorithm is actually more closely related to ``coupling into
and from the past.''  We say more about this extension of CFTP and these
two algorithms in \sref{ciaftp}.

The next development was ``monotone-CFTP''
\citep{propp-wilson:exact-sampling}, which is a particularly efficient
algorithm that can be used when the state space has a partial order
$\preceq$ that is preserved by the randomizing operations (if
$x\preceq y$ then $f_{-t}(x)\preceq f_{-t}(y)$), and there is a
biggest state $\hat{1}$ and smallest state $\hat{0}$.
These conditions are somewhat restrictive, but a
surprisingly wide variety of Markov chains of practical interest
satisfy these conditions; see e.g.\ the examples given by
\cite{propp-wilson:exact-sampling},
\cite*{luby-randall-sinclair:markov-lattice},
\cite{felsner-wernisch:linear-extension},
\cite*{haggstrom-lieshout-moller:exact-spatial},
\cite{lund-wilson:storage},
\cite{berg-steif:codings},
\cite{nelander:beach},
\cite*{mira-moller-roberts:slice}, and 
\cite*{muri-chauveau-cellier:dna}.
The algorithm in \fref{citp} computes compositions
in the order
$$(\cdots((f_{-1} \circ f_{-2}) \circ f_{-3}) \cdots \circ f_{-T+1}) \circ f_{-T}.$$
For monotone-CFTP (and most subsequent versions of CFTP),
it is much easier to perform the composition in the order
$$f_{-1} \circ (f_{-2} \circ (f_{-3} \cdots \circ (f_{-T+1} \circ f_{-T}) \cdots )).$$
The reason is that in the end we only need the image of the final
composition, and if we compose the maps in the second order, then we
only need to compute the images of the intermediate compositions,
rather than having to compute the entire map.  (We explain below how
these images are computed --- what's important here is that this computation
is easy in the monotone setting.)  In contrast with the first order
of compositions, where we compose maps going back in time, doing
compositions in the second order requires us to pick a starting value
$-T$ in some fashion, and compose maps going forwards in time back to
the present.  If the composition is not coalescent (i.e.\ the image is
not a singleton), then we pick another starting value even further
back in the past.  A reasonable choice of starting times are times of
the form $-2^k$, and the resulting binary-backoff version of CFTP is
shown in \fref{binary-backoff-cftp}.

Note that the algorithm resets its source of randomness in a manner
that ensures that for each $t$, the random map $f_{-t}$ has the same
value each time it is used in a composition.  {\it A priori\/} we should
be extremely suspicious of any proposal to pick fresh values for $f_{-t}$
each time it is refered to, since then the binary-backoff CFTP in
\fref{binary-backoff-cftp} would not properly emulate the algorithm in
\fref{citp}.  In fact it is a bad idea, and results in a biased
algorithm.  For this reason it is emphasized that the same coins need
to be re-used each time that $f_{-t}$ is generated and refered to, and
it becomes unclear how to proceed with a read-once source of
randomness.

\begin{figure}[htb]
\begin{code}
\item   BinaryBackoffCFTP (NumberOfSamples) 
\item\s   \For $i := 1$ \To NumberOfSamples \{
\item\s\s   $T := 1$
\item\s\s   \Repeat \{
\item\s\s\s   Set := \sspace
\item\s\s\s   \For $t$ := $T$ \DownTo $1$
\item\s\s\s\s   \If $t$ is a power of 2
\item\s\s\s\s\s   SetRandomSeed(seed[$i$,$\log_2(t)$])
\item\s\s\s\s   ApplyRandomMap(Set)
\item\s\s\s   $T := 2*T$
\item\s\s   \} \Until Singleton(Set)
\item\s\s   \Output ElementContainedIn(Set)
\end{code}
\caption{
  Pseudocode for the binary-backoff implementation of CFTP, iterated
  some number of times.  The variable Set represents the image of the
  composition $f_{-T}\circ\cdots\circ f_{-t}$, or more generally a
  superset of the image.  The representation of Set is seldom a naive listing
  of states.  Note that a number of independent random
  number streams are used, and some streams of randomness are read
  multiple times.  }
\label{fig:binary-backoff-cftp}
\end{figure}

{\it Remark:\/}  The spanning tree algorithm and the dead-leaves
process are unusual in that they compose their maps using the first
order, i.e.\ back into the past rather than from the past.  As pointed
out by Kendall, these algorithms therefore
already run with a read-once source of
randomness.  Nearly every other CFTP-type algorithm composes maps
forward in time from the past, and therefore requires a different
method of running with read-once source of randomness.

We briefly return to monotone-CFTP and explain how it computes the
images of the compositions of random maps when they are composed going
forwards in time.  Technically the
precise image of the map is not computed, but rather a superset of the
image is computed.  The superset at time $-t$ is represented in the
computer by two bounding states, $\ell_{-t}$ and $u_{-t}$, and the
superset is the interval $\{x:\ell_{-t} \preceq x \preceq u_{-t} \}$.  The
bounding values $\ell_{-T}$ and $u_{-T}$ are set to the
minimum and maximum states respectively, so that trivially the resulting
interval is a superset of the image of the initial composite map
(indeed of any map).  The bounding states are updated by the rules
$u_{-t+1} := f_{-t}(u_{-t})$ and $\ell_{-t+1} := f_{-t}(\ell_{-t})$.  Since
the random maps respect the partial order of the state space, by
induction we see that the interval $\{x:\ell_{-t+1} \preceq x \preceq u_{-t+1}
\}$ must be a superset of the image of the map $f_{-T}\circ\cdots\circ f_{-t}$.
We remark that even though the image does not necessarily occupy the
entire interval, the image is a singleton if and only if the interval
is a singleton (i.e.\ iff $\ell_t=u_t$).

After the success of monotone-CFTP, there has been a good deal of
research on finding more classes of applications to which CFTP may be
efficiently applied; see e.g.\
\cite{propp-wilson:unknown-markov-tree},
\cite{kendall:area-interaction},
\cite{kendall:boolean},
\cite{haggstrom-nelander:antimonotone},
\cite{luby-vigoda:independent},
\cite{murdoch-green:continuous},
\cite{kendall-moller:exact-spatial},
\cite{moller:conditional},
\cite{haggstrom-nelander:random-fields},
\cite{green-murdoch:bayesian},
\cite{kendall-thonnes:geometry},
\cite{huber:techniques},
\cite{huber:swendsen-wang}, and
\cite{haggstrom-steif:finitary}.
In this subsequent work, researchers have studied state spaces that
don't have a convenient partial order preserved by the random maps,
and found other clever mechanisms for effectively representing and
updating a superset of the image of the composition of the random
maps.  There is a tradeoff in the choice of representation:
maintaining the exact image or a very detailed superset of it may take
a lot of computer effort, while if too course a superset is
maintained, coalescence may not be readily detected.

We remark that there are also a number of perfect sampling algorithms
based on ``Fill's algorithm'' \citep{fill:interruptible} rather than
CFTP (see e.g.\
\cite{fill:interruptible},
\cite{thonnes:points},
\cite{moller-schladitz:fill}, and
\cite*{fill-machida-murdoch-rosenthal:fill}), and that there are Markov chain-based
perfect sampling algorithms based on neither method (see e.g.\ 
\cite*{asmussen-glynn-thorisson:unknown-markov},
\cite{aldous:unknown-markov},
\cite{lovasz-winkler:unknown-markov}, and
\cite{propp-wilson:unknown-markov-tree}).

\old{
The procedure attempts to determine the value of After the while loop
has executed $t$ times, the variable $F$ represents the random mapping
defined by the composition of the last $t$

which is sometimes implementable more or less directly (e.g.\ the 
dead-leaves process), but which is more often implemented by
Some versions of CFTP already work with a read-once stream of random coins,
such as the dead-leaves process, which is nicely animated at
 {\tt http://www.warwick.ac.uk/\char126statsdept/Staff/WSK/dead.html}
I don't mean that just *some* versions of CFTP can be done with a read-once
stream of randomness, what I mean is

Here the variable $S$ is a computer representation of (a superset of)
the image of $F_T^t$, which is the composition of the random maps from
time $T$ (inclusive) to time $t$ (exclusive).  ``\sspace'' is almost
never a naive listing of all possible states.  In the case of
monotone-CFTP, the representation of the image of $F_T^t$ consists of
just two states, and the image is a singleton set precisely when these
two states are equal.  Other versions of CFTP use more involved
schemes to represent a not-too-large superset of this image in a
reasonably concise manner.  ApplyRandomMap(S) generates a new random 
map using the current stream of randomness, and updates the representation
of the image set $S$ appropriately.
}

\section{Read-once CFTP}
\label{sec:rocftp}

In this section we explain the read-once randomness version of CFTP,
for which pseudocode is given in \fref{read-once-cftp}.  Read-once
CFTP may be viewed as a retroactive stopping rule.  It applies random
maps going forwards in time, and then at some point it decides to
stop, and then returns not the current state, but some previous state.

A key part of read-once CFTP is a composite random map procedure,
which uses the ApplyRandomMap procedure as a subroutine.  From the
standpoint of read-once CFTP, it appears as if the composite map
procedure generates a random map, makes some effort to determine
whether or not the map is coalescent (i.e.\ whether or not it maps
all states to one
state), and then evaluates the map at a given input state to obtain an
output state.  The composite random map preserves the desired
probability distribution, in that if the input state is distributed
according to the desired distribution, then so is the output state.
If the procedure determines (by examining the representation of the
superset of the image of the map) that the random map is coalescent,
then we say that the map is ``officially coalescent.''  Otherwise the
map is not officially coalescent, and it may or may not map all states
to one state.  It is important that the choice of input state at which
the random map is evaluated does not affect whether or not the
composite map procedure detects coalescence (since otherwise it would not
appear as if the procedure tested for colescence and then evaluated
the random map at the input state).  For efficiency reasons,
we design the procedure so that it produces an officially coalescent
random map with probability $p\geq 1/2$.  We assume that subsequent
invocations of the composite map procedure are independent.  Later we
explain how to implement a composite random map that meets these
requirements, but first we see how to use it for read-once CFTP.

Suppose the composite map procedure gave us the entire random map,
rather than just evaluating it at one state.  Then
we could do CFTP, composing new composite
maps going back in time.  Let $f_{-T}$ be the first (closest to time $0$,
smallest $T$) composite map that is officially coalescent.  $T$ is a
geometric random variable with mean $1/p\leq 2$.  $f_{-T}$ is a random
composite map conditioned to be officially coalescent, and is
furthermore independent of $T$.  Let $S$ be the state in the image of
$f_{-T}$.  CFTP would then apply the composite maps $f_{-T+1}, \ldots,
f_{-1}$ to $S$, and return the result.  The composite maps $f_{-T+1},
\ldots, f_{-1}$ are i.i.d.\ random composite maps conditioned not to
be officially coalescent, and are independent of $S$.  So we could
equivalently generate $T-1$ fresh random composite maps conditioned
not to be officially coallescent, and apply them to $S$.  Furthermore,
there is no need to count $T$.  We can simply update $S$ using fresh
composite random maps, until one of the maps is officially coalescent,
and return the value of $S$ prior to the last composite map.  (This is
where we use the i.i.d.\ condition.)

Thus to generate a random sample, we make random composite maps until
we see two that are officially coalescent, and compose those maps
between the first coalescent map (inclusive) and the second coalescent
map (exclusive).  Since the second of the officially coalescent
composite maps is used only as a stopping criterion, and is not itself
included in the composition of maps which results in the random
sample, this second coalescent map is independent of the returned
random sample, and so may be used in the generation of a subsequent
independent random sample.  If $k$ random samples are desired, the
last sample is returned upon the generation of the $(k+1)^{\text{st}}$
officially coalescent composite map.

\begin{figure}[htb]
\begin{code}
\item   ReadOnceCFTP (NumberOfSamples)
\item\s   Initialize()
\item\s   \For $i := 1$ \To NumberOfSamples
\item\s\s   \Output NextSample()
\item
\item Initialize ()
\item\s   State := $\langle\text{arbitrary state}\rangle$
\item\s   \Repeat
\item\s\s   ApplyCompositeMap(State,CoalescenceFlag)
\item\s   \Until CoalescenceFlag
\item
\item NextSample ()
\item\s   \Repeat
\item\s\s   OldState := State
\item\s\s   ApplyCompositeMap(State,CoalescenceFlag)
\item\s   \Until CoalescenceFlag
\item\s   \Return OldState
\end{code}
\caption{
Pseudocode for the read-once version of CFTP, iterated some number of
times.  The Initialize() routine generates a random map conditioned to
be officially coalescent.  The NextSample() routine uses this map,
and produces a random
sample from the desired distribution, as well as an independent
random map conditioned to be officially coalescent, to be used
in the subsequent call to NextSample().  The ApplyCompositeMap()
procedure overwrites the input state
with the output state, and stores in CoalescenceFlag a Boolean value
reporting whether or not the random composite map is officially
coalescent.  This version of CFTP uses only one pseudorandom stream, and
random values never need to be re-read.}
\label{fig:read-once-cftp}
\end{figure}

\break

If there were some standard notation for reasoning about algorithms
that produce random outputs, we might be able to re-express the
previous discussion more symbolically in a manner such as the
following.  Here we have assumed that there is a positive probability
that a random (composite) map is coalescent, and we have let $f(*)$
denote the unique element in the image of a coalescent map $f$.

\begin{align*}
\left\{\begin{tabular}{l}
 $X :=$ draw from $\pi$\\
 $f :=$ officially coalescent random map\\
 \Output $X,f(*)$
\end{tabular} \right\}
&\stackrel{\cal D}{=}
\left\{\begin{tabular}{l}
 $T:=0$\\
 $F := \langle\text{identity map}\rangle$\\
  \While $F$ is not coalescent\\
  \s $T := T+1$ \\
  \s $f_{-T} := \rmap$\\
  \s $F := F \circ f_{-T}$ \\
 $f :=$ officially coalescent random map\\
 \Output $F(*),f(*)$
\end{tabular} \right\} \\
&\stackrel{\cal D}{=}
\left\{\begin{tabular}{l}
 $T:=0$\\
 \Repeat\\
  \s $T := T+1$ \\
  \s $f_{-T} := \rmap$\\
 \Until $f_{-T}$ is officially coalescent\\
 $f :=$ officially coalescent random map\\
 \Output $f_{-1}(\cdots f_{-T+1}(f_{-T}(*))),f(*)$
\end{tabular} \right\} \\
&\stackrel{\cal D}{=}
\left\{\begin{tabular}{l}
 $T:=0$\\
 \Repeat\\
  \s $T := T+1$ \\
  \s $f_{-T} := \rmap$\\
 \Until $f_{-T}$ is officially coalescent\\
 $f :=$ officially coalescent random map\\
 \Output $f_{-T+1}(\cdots f_{-1}(f(*))),f_{-T}(*)$
\end{tabular} \right\} \\
&\stackrel{\cal D}{=}
\left\{\begin{tabular}{l}
 $f :=$ officially coalescent random map\\
 State $:= f(*)$ \\
 \Repeat\\
  \s OldState := State \\
  \s $f := \rmap$\\
  \s State := $f($State: Exp $)$\\
 \Until $f$ is officially coalescent\\
 \Output OldState, State
\end{tabular} \right\}
\end{align*}

In the\old{ implementation of the} composite map procedure given in
\fref{interleaved}, we independently and in parallel
update two subsets of the state space, each representing (a superset
of) the image of a random map.  Initially the two maps are the
identity map, so that the two subsets are initially the whole state
space.  At each step we update the first set with ApplyRandomMap, and
update the second set similarly but with an independent random map.
We keep doing these parallel updates until the second map is officially
coalescent.  The number of times that the first subset was updated is
independent of the mappings used to do its updates.  Therefore the
first mapping preserves the desired probability distribution on the
state space.  Furthermore, since it is with probability at least $1/2$
that the first subset becomes a singleton no later than the second
subset, it is with probability at least $1/2$ that the first mapping
is officially coalescent.  (This is where we used condition 2, coalescence
in finite time.)

\begin{figure}[htb]
\begin{code}
\item   ApplyCompositeMap (State,CoalescenceFlag) 
\item\s   Set1 := \sspace
\item\s   Set2 := \sspace
\item\s   \While \Not Singleton(Set2) 
\item\s\s   ApplyRandomMap(Set1,State) /* apply same random map to Set1 and State  */
\item\s\s   ApplyRandomMap(Set2) \ \ \ \ \ \ /* but apply independent random map to Set2 */
\item\s   CoalescenceFlag := Singleton(Set1)
\end{code}
\caption{
Interleaved version of ApplyCompositeMap().
ApplyRandomMap() takes an optional argument State
that is updated according to the generated random map.} 
\label{fig:interleaved}
\end{figure}

\section{Read-once CFTP and PASTA}
\label{sec:pasta}

We obtained this read-once version of CFTP by starting with the usual
version of CFTP and modifying it.  It would also have been possible to
start with what is known as PASTA in the queuing theory and operations
research literature, make suitable changes, and arrive at read-once
CFTP.  In this section we explain PASTA and this alternate derivation of
read-once CFTP.

PASTA is a statement about a stochastic process evolving in time, and
discrete events which affect the stochastic process and occur at times
given by a Poisson process.  PASTA stands for ``Poisson arrivals see
time averages,'' which means that the steady-state distribution of the
stochastic process averaged over all times is identical to the
steady-state distribution of the process sampled at the times just
prior to the Poisson events.  \cite{wolff:pasta} introduced the
concept of PASTA, and showed that it holds whenever the stochastic
process cannot anticipate the future driving events.  Since that time
there have been many articles on applications and generalizations of
PASTA, which go by a variety of acronyms, including ASTA, ESTA, EATA,
EPSTA, CEPSTA, and MUSTA; reviews are given by
\cite{melamed-whitt:pasta}, \cite*{bremaud-kannurpatti-mazumdar:pasta},
and \cite{melamed-yao:pasta}.

A discrete-time version of PASTA would be a statement about a discrete
time Markov chain, and random events that occur at integer times.  If
there is an event at a given time, then the next state of the Markov
chain is drawn according to one transition rule, while if there is no
event, then a different transition rule is used.  Discrete-PASTA would
state that the distribution of the Markov chain sampled at times just
prior to when events occur will be identical to the steady-state
distribution of the Markov chain.

In read-once CFTP, an event occurs precisely when a composite map is
officially coalescent.  Imagine first randomly picking those integers
at which events occur.  If there is an event at a given time, then the
Markov chain is updated by a random composite map conditioned to be
officially coalescent, otherwise it is updated by random composite map
conditioned not to be officially coalescent.  Discrete-PASTA asserts
that the if we draw samples from the Markov chain at times just prior
to when the composite maps are officially coalescent, the steady-state
distribution of the draws will be the steady-state distribution of the
Markov chain.  PASTA is a statement about the steady-state behavior of
the draws; in general the first several draws taken at positive times
will be out of equilibrium.  In this particular application of PASTA,
since there is a coalescent map between draws, not only are draws
after the first one easy to compute, but they also must necessarily be
independent of one another.  Since the draws are independent, any
particular draw is already in the steady-state distribution.
Read-once CFTP ignores the first draw (since it is neither in
equilibrium nor easy to compute), and outputs the subsequent draws
until the desired number of independent perfectly random samples are
generated.

We remark that CFTP and PASTA are not completely unrelated ideas.  The
``time zero sees time averages'' principle behind CFTP can be used to
derive the ``Poisson arrivals see time averages'' in both the continuous
and discrete settings.  Perhaps further connections can be made
between perfect simulation algorithms and the various generalizations
of PASTA.

\section{Performance of read-once CFTP}
\label{sec:performance}

\subsection*{Expected running time}

Let $T_N$ be the expected number of random maps that we need to compose
before the composition is officially coalescent.  We should not hope to
have an algorithm that uses many fewer than $T_N$ random maps.
The usual binary-backoff CFTP
uses between $2 T_N$ and $4 T_N$ random maps, with the constant typically
being around $2/\log 2 \approx 2.9$.  The read-once version of CFTP will use on
average $4(k+1) T_N$ random maps to generate $k$ samples: The expected
time to generate a composite map is $2 T_N$.  The
time to generate an officially coalescent composite map is $\leq 4 T_N$,
and we do this $k+1$ times to generate $k$ samples.
%
If updating a single state is about
as expensive as updating a whole set of states, then since each
composite map does on average $T_N$ 
state updates in addition to the expected $2 T_N$ set updates,
we expect to do $6 (k + 1) T_N$ updates altogether.

In some applications of CFTP, particularly on continuous state spaces,
applying the first several random maps can be enormously more
expensive than applying subsequent random maps, because initially the
updated set is the whole state space, and later it is smaller.
The binary-backoff version of CFTP does these expensive updates a
number of times that is logarithmic in $T_N$, while the read-once
version of CFTP does the expensive updates on average $4$ times (plus
another $4$ times for the first sample).  Therefore, for these
applications, the read-once version of CFTP could be up to
logarithmically faster than the binary-backoff version of CFTP.  This
issue of logarithmic slowdown when the initial updates are expensive
has come up before.  For instance, one algorithm given by
\cite{propp-wilson:unknown-markov-tree} used random maps with this
sort of run time variability, and used composite maps to avoid the
logarithmic slowdown.

\subsection*{Tail distribution of the running time}

Suppose $S$ and $T$ are subsets of the state space, and $f$ denotes a
random map.  It is clear that if $S\subseteq T$ then $f(S)\subseteq
f(T)$.  The ApplyRandomMap procedure when applied to $S$ updates $S$
to a superset of $f(S)$, so conceivably this superset of $f(S)$ may
not be contained within the corresponding superset of $f(T)$.  But in
practice every ApplyRandomMap procedure that anyone uses respects
subset inclusion, meaning that the superset of $f(S)$ is contained
within the superset of $f(T)$.  \cite{kendall-moller:exact-spatial}
call this property the funnelling property.

Under the assumption that the ApplyRandomMap procedure satisfies the
funnelling property, the tail distribution of the running time decays
geometrically with decay constant that is a (universal) constant
multiple of $T_N$.  A similar upper bound holds for the usual CFTP,
but if the underlying Markov chain has a sharp threshold, the tail
distribution could be even tighter.

As mentioned above, the new version is ``temptation free,'' whereas
the usual version occasionally enters states where the user may be
tempted to abort and restart (provided that the funnelling property
holds in that the random maps each take the same amount of time to
apply).  The temptation-free property holds provided that the user
does not look at the value of the counter, or if the user might do
such a thing, the alternative random map procedure in
\fref{interleaved} can be used instead since it has no counter.
Because of the funnelling property, the number of iterations before
the composite map procedure returns is always stochastically dominated
by the number of iterations required by a fresh call to
ApplyCompositeMap.  Furthermore, the number of calls to apply
composite map before the next several samples are returned is
stochastically dominated by the number of such calls if the user were
to restart ReadOnceCFTP.  Therefore, under these assumptions about
ApplyRandomMap, the user will never get his or her desired samples
more rapidly by interrupting and restarting the ReadOnceCFTP
procedure.  (As mentioned by \cite{propp-wilson:unknown-markov-tree}
and \cite{fill:interruptible}, for some applications the underlying
random maps take a variable amount time to apply.  For these
applications one should not expect ReadOnceCFTP to yield a
temptation-free sampling algorithm, nor should one expect Fill's
algorithm to yield an interruptible sampling algorithm.)

\subsection*{Memory}

The memory required for the binary-backoff version of CFTP is the
memory to store a subset, plus the memory to store two integers ($T$
and $t$).  The memory required for read-once CFTP 
using the interleaved
version of the composite random map procedure is twice the
memory to store a subset plus the memory to store a state.  In
principle there is no upper bound on the storage requirements for
these integers, so in principle the new version of CFTP could have
smaller memory requirements, but in practice finding memory for these
integers is a non-issue.

More significant is the effect of the constant factor increase in
memory requirements associated with storing two subsets of the state
space.  Computers typically contain several different types of memory,
including an L1 cache, an L2 cache, and a main memory composed of
DRAM.  The memory close to the processor is fast, expensive, and
small, while the main memory is slow, cheap, and large (see e.g.\ 
\cite[Chapter 7]{hennessy-patterson:architecture}).  Even if the
simulation still fits within main memory, if less of it fits within
the caches, performance will degrade.  \cite{yan:memory} did timing
experiments of a wide variety of sizes of Ising model simulations, and
reported that it was quite noticeable when the next slower type of
memory started to be used.

Therefore, unless the memory requirements are quite small, we
recommend instead the version of the composite map procedure given in
\fref{composite}.  Rather than updating the two sets in parallel, only
one set is updated, and then later only the other set is updated.
Only one subset, a state, and an integer need to be stored.
This version of the procedure behaves in the same manner as the
interleaved version, unless the counter overflows.  Even if the count
were to overflow the integer, while the run time performance could be
affected slightly, the distribution of the output of the algorithm is still
identical to the desired distribution.

\begin{figure}[htb]
\begin{code}
\item   ApplyCompositeMap (State,CoalescenceFlag) 
\item\s   Set := \sspace
\item\s   Count := 0
\item\s   \While \Not Singleton(Set) 
\item\s\s   ApplyRandomMap(Set)
\item\s\s   Count := Count$+1$
\item\s   Set := \sspace
\item\s   \While Count$>$0 
\item\s\s   ApplyRandomMap(Set,State) /* apply same random map to Set and State  */
\item\s\s   Count := Count$-1$
\item\s   CoalescenceFlag := Singleton(Set)
\end{code}
\caption{
Memory-efficient version of ApplyCompositeMap().}
\label{fig:composite}
\end{figure}

\subsection*{Overall}

In some circumstances, but certainly not all, it may be preferable
to use read-once CFTP.

\section{Read-once CFTP and unbounded state spaces}
\label{sec:unbounded}

In this section we describe a small modification to read-once CFTP
that makes it easier to use with unbounded state spaces.

For some applications of CFTP to sampling from unbounded state spaces,
it is convenient to mix two different Markov chains on the same state
space.  For instance, \cite{murdoch:exact-bayesian} describes examples where the
natural Markov chain for a state space has favorable mixing properties when started from most typical states, but that when
started from points ``very far away'' in the tails of the stationary
distribution, the time to randomize can get arbitrarily large.  Such a
Markov chain is ``non-uniformly ergodic,'' and it has been observed by
a number of authors (for example
\cite{foss-tweedie:srs-cftp}) 
that if we do CFTP using such a Markov chain in a
straightforward fashion, coalescence takes infinitely long.
(The reason is that the coupling time
upper bounds the worst case mixing time (see e.g.\
\cite{aldous:group-walks}), which is infinite for non-uniformly ergodic
Markov chains.)  

To speed
up the convergence time to a finite value, \cite{murdoch:exact-bayesian}
suggested mixing the natural Markov chain with another Markov chain
called the ``independence sampler.''  Details of how to do this can be
found in \citep{murdoch:exact-bayesian}; the algorithm for point processes
in \sref{point-process} also serves as an illustrative example.
We mention here that the effect of the
independence sampler is to map the entire state space to a bounded
region, but otherwise the independence sampler has poor convergence
properties.  If an algorithm occasionally makes moves using the
independence sampler, but most of the time using the natural Markov
chain, then the convergence time will be finite, and reasonably fast
for the examples considered by Murdoch.  Murdoch's solution is fairly effective, and upon learning of it, \cite{wilson:multishift} used it in a perfect sampling algorithm for the autonormal distribution.

\cite{murdoch:exact-bayesian} originally suggested flipping a suitably
biased coin at each time step to decide whether to update using the
independence sampler or the natural Markov chain.  But for the point
process example in \sref{point-process}, if the independence sampler
is applied to frequently, it tends to disrupt coalescence detection.
For the autonormal algorithm given by \cite{wilson:multishift},
independence sampler updates are more expensive than normal updates.
But if the independence sampler is used too infrequently, the expected
run time is guaranteed to be large.  An alternative is to let the
composite map procedure determine on its own what the right mixing
ratio is.

Our recommendation for applications using the independence sampler is to let
read-once CFTP's composite random map procedure do one
update from the independence sampler, and do subsequent updates
using the natural Markov
chain.  This change is most easily made by replacing the lines which
initialize Set to the whole state space with lines that
instead initialize it to the result of the first random map,
as shown in \fref{rocftp-1st}.  For both the autonormal and point process
applications, using the independence sampler for only the first update
also helps simplify the code.

\begin{figure}[htb]
\begin{code}
\item   ApplyCompositeMap (State,CoalescenceFlag) 
\item\s   Set := ImageOfFirstRandomMap()
\item\s   Count := 0
\item\s   \While \Not Singleton(Set) 
\item\s\s   ApplyRandomMap(Set)
\item\s\s   Count := Count$+1$
\item\s   Set := ImageOfFirstRandomMap(State)
\item\s   \While Count$>$0 
\item\s\s   ApplyRandomMap(Set,State) /* apply same random map to Set and State  */
\item\s\s   Count := Count$-1$
\item\s   CoalescenceFlag := Singleton(Set)
\end{code}
\caption{
Version of composite map procedure suitable for many applications
of read-once CFTP on unbounded
state spaces.  The first random map may implement a fundamentally
different Markov chain with the same stationary distribution, or in
some cases the first map may be distributed in the same manner as
subsequent random maps, but it is convenient to treat it differently.
Like ApplyRandomMap(), ImageOfFirstRandomMap() takes an optional
argument State that is updated according to the generated random map.
} 
\label{fig:rocftp-1st}
\end{figure}

{\it Remark:\/}
Since the random maps within the composite map procedure are no longer
all identically distributed, it is no longer {\it automatic\/} that the tail
distribution of the running time decays exponentially.  For the
applications to the autonormal \citep{wilson:multishift} and to point
processes (\sref{point-process}), and
perhaps for other applications, it is elementary to show that the
tails still decay exponentially.  But conceivably for some application
the distribution could have fat tails, and the expected running time
could even be infinite.  Under such conditions,
\cite*{luby-sinclair-zuckerman:speedup} recommend restarting (with the
independence sampler) after runs (of the natural Markov chain) of
lengths $1,1,2,1,1,2,4,1,1,2,1,1,2,4,8,1,\dots$.  They showed that this
sequence is a ``universal restart sequence'' that can be used to
speed up the expected time for an event to happen when the running time
distribution has fat tails.  But until fat-tailed coalescence time
distributions show up in perfect simulation, we recommend that only
the first update be made using the independence sampler.

In other applications of CFTP on unbounded state spaces, it is
convenient to implement the first random map in a different manner
than subsequent random maps, even though from a mathematical
standpoint the random maps themselves are drawn from the same
distribution.  For instance, \cite*{haggstrom-lieshout-moller:exact-spatial}
consider the Widom-Rowlinson model on a finite region (such as a unit
square) together with a monotone Markov chain.  With probability 1, a
random state will consist of finite number of red points and a finite
number of blue points from this region.  In the natural partial order,
we have $X\preceq Y$ if each red point in configuration $X$ is also a
red point of configuration $Y$, and each blue point of configuration
$Y$ is also a blue point of configuration $X$.  There is no top state
or bottom state in this partial order if we restrict our attention to
sets with finitely many red and blue points, which is inconvenient if
we wish to run monotone-CFTP.  But there are top and bottom states
with very geometric interpretations: in the top state each point in
the region is red and no points are blue, and vice versa for the
bottom state.  Conveniently, after applying one random map, each state
in the image of the map is sandwiched in between an upper state and a
lower state, both of which with probability 1 contain only finitely
many points.  From an implementation standpoint, it is inconvenient
(though possible) to create a data structure that can represent finite
collections of red and blue points, as well as the uncountable
collections of points that we would need to represent the top state in
bottom state.  Since after one random map we only need to represent
finite collections of points, it is sensible to choose a simpler data
structure that can only represent finite collections of points, and
treat the first random map as a special case.  Another application for
which it is convenient (though not necessary) to treat the first
random map as a special case, on account of the state space being very
large, is the autogamma sampler given by \cite{moller:conditional}.

To run read-once CFTP on applications for which it is convenient to
treat the first random map has a special case, as before
(\fref{rocftp-1st}), we replace the lines initializing Set
to the whole state space with lines initializing it to the
result of the first random map.  Since the random maps within
the composite map are still i.i.d.\ even though they are implemented
differently, it is once again automatic that the tail distribution
of the running time decays geometrically.

\section{Background on coupling into and from the past}
\label{sec:ciaftp}

The CFTP-type algorithms which do not use independent random maps
(i.e.\ don't satisfy condition 1 of the main result) are the
``coupling {\it into and from\/} the past'' (CIAFTP) algorithms.  These
algorithms extend CFTP, and were introduced by
\cite{kendall:area-interaction}, though he did not use this term.
There are comparitively few applications of coupling into and from the
past, as opposed to coupling {\it from\/} the past, but new ones may
be developed as more people become aware of this worthwhile technique.
In the next two sections we adapt CIAFTP algorithms to the setting of
read-once randomness, and as preparation, we review the basic method
here.

For concreteness, we explain coupling into and from the past by means
of an example, which we then generalize and modify.  Recall that in
\sref{cftp} we mentioned that the dead-leaves process and the spanning
tree algorithm of \cite{aldous:tree} and \cite{broder:tree} were both
based upon the ``state at time zero is random'' principle.  In the
same way that (as \cite{kendall-thonnes:geometry} point out)
the dead leaves process can
be regarded as an early form of CFTP, where the random maps are
composed going backwards in time into the past rather than forwards in
time from the past, the Aldous/Broder spanning tree algorithm is an
early form of coupling into and from the past.  The comparison
``dead-leaves : CFTP :: spanning tree : CIAFTP'' is sufficiently
compelling that we explain both algorithms together.

As Broder and Aldous explain in their writeups of the spanning tree
algorithm, there are two different Markov chains that are run together
in a coupled fashion.  The target Markov chain (that we wish to sample
from) is on the set of rooted spanning trees of a given graph.  The
other Markov chain, which we shall call the reference chain, is the
simple random walk on the given graph.  It is assumed that we already
know how to sample from the reference chain.  In the applications of
coupling into and from the past given by
\cite{kendall:area-interaction}, \cite{kendall-moller:exact-spatial}, and
\cite{lund-wilson:storage}, the chain that we already know how to
sample from is called the dominating chain, since its values
stochastically dominate the values of the target chain.  In the
spanning tree application however, there is no natural partial order,
or at least none that anyone has found, so the term ``dominating
chain'' is not appropriate in general.

The target chain for rooted spanning trees moves the root to a random
neighboring vertex, adjoins an edge directed from the old root to the
new root, and then removes the edge directed out of the new root.  It
is an interesting exercise to verify that this Markov chain preserves
the uniform distribution on rooted spanning trees.

The coupling between the target chain and the reference chain is such
that the random walk on the graph follows the same trajectory as the
root of the spanning tree.  The algorithm picks a random value for the
reference chain, runs it backwards in time, and attempts to determine
the state of the target chain at time $0$.  Since the algorithm can
use the information provided by the reference chain, it is only
necessary to test if those spanning trees with a specified root,
rather than all rooted spanning trees, coalesce to a single spanning
tree by time $0$.  The coupling used to determine this is actually
quite similar to the coupling used in the dead leaves process.  We can
view the rooted spanning tree as a vector assigning each vertex to its
parent, and the dead leaves process as a vector (indexed by $\R^2$!)
describing how the dead leaves tesselation looks at each point.  As
the target chain runs forward, parts of the vector are overwritten
with new values, and the remaining parts are left untouched.  For the
spanning trees, the vector is changed only at the old root and the new
root.  This overwriting process is inherently Markovian, and is
determined by the reference chain.  For the dead leaves process, the
overwriting is determined by an i.i.d.\ process.  It is easy to
directly compose backwards in time the random maps determined by an
overwriting process: initialize with a clean slate, and then only
overwrite those portions that have not yet been touched.  Keep going
until every part of the vector has been touched.

Since simple random walk on an undirected graph is reversible, it is
easy to run the reference chain backwards in time and use it to
determine what the overwriting process did in the past.  Later
\cite*{kandel-matias-unger-winkler:shuffling-biological} had reason to
generate random spanning trees from an Eulerian directed graph, i.e.\
a directed graph where the in-degree of any vertex is also its
out-degree.  While simple random walk on the Eulerian graph is no
longer reversible, they pointed out that the time reversal of this
directed random walk is still easy to simulate, so that essentially
the same method can be used to generate random spanning trees on
directed Eulerian graphs.  (Other tree algorithms that work for more
general directed graphs are given by Propp-Wilson.)

If we abstract away the particulars of the spanning tree algorithm
while maintaining the overall strategy, we get the ``coupling into and
into the past'' procedure, for which pseudocode is given in
\fref{ciaitp}.
The algorithm generates a sequence of states
$\ldots,X_{-3},X_{-2},X_{-1},X_0$
of the reference Markov chain
and a sequence of random maps
$\ldots,f_{-3},f_{-2},f_{-1}$ of the target Markov chain, so that for
each time $-T$ in the past, the following two properties hold
\begin{enumerate}
\item $X_{-T}$ is distributed according to $\pi_{\text{ref}}$.
\item The pairs $(X_{-t+1},f_{-t})$ (for $-T\leq-t\leq-1$) look as
if they were generated by the ``useful coupling'' between the reference
Markov chain and random maps of the target Markov chain.
\end{enumerate}
Naturally, for any given value of $-T$, it would be straightforward to
simply start at time $-T$ and generate the pairs $(X_{-t+1},f_{-t})$
running forwards in time.  As with the usual version of CFTP, the
algorithm attempts to determine the state of the target Markov chain
at time $0$.  If there is only one possible value for the state at
time $0$, then this state is a draw from the stationary distribution
of the target chain.  Of course the algorithm does not know ahead of
time what time $-T$ to start with, so it must be able to generate  
the states $\ldots,X_{-3},X_{-2},X_{-1},X_0$ of the reference chain
and the random maps $\ldots,f_{-3},f_{-2},f_{-1}$ of the target chain
going backwards in time rather than forwards.

\begin{figure}[phtb]
\begin{code}
\item   $X :=$ ReferenceChainRandomState()
\item   $F := \langle\text{identity map on target chain}\rangle$
\item   \While \Not Singleton(ImageOf($F$ restricted to states compatible with $X$))
\item\s   $X' :=$ ReverseReferenceChain($X$)
\item\s   $F := F \circ$ TargetChainRandomMapCoupledExPostFacto($X'$,$X$)
\item\s   $X := X'$
\item   \Return ElementContainedIn(ImageOf($F$ restricted to states compatible with $X$))
\end{code}
\caption{
High level pseudocode for coupling into and from the past.  The state
$X$ of the reference chain may be used when determining the image of
the map $F$.  Since both the target chain random maps are composed
going backwards in time, in addition to the reference chain being run
backwards, the algorithm might be more properly called coupling {\it
into and into\/} the past.}
\label{fig:ciaitp}
\end{figure}

To generate the states $X_{-t}$ and random maps $f_{-t}$ going
backwards in time, first the state $X_0$ is drawn from the stationary
distribution $\pi_{\text{ref}}$ of the reference Markov chain.  Then the
time-reversal of the reference Markov chain is run, thereby producing the
requisite sample path $\ldots,X_{-3},X_{-2},X_{-1},X_0$ of the reference
chain up to time $0$.  (The reference chain is typically
reversible, so that running its time-reversal is the same as running
the original Markov chain.)  Implicit in property 2 above is that the
conditional distribution of $f_{-t}$ is a function of $X_{-t}$ and
$X_{-t+1}$ alone.  Recall that given the state $X_{-t}$ of the reference
chain at time $-t$, the pair $(X_{-t+1},f_{-t})$ is supposed to be
distributed according to the presepecified ``useful coupling''.  Since
both $X_{-t}$ and $X_{-t+1}$ are generated before the random map
$f_{-t}$, the map $f_{-t}$ must be coupled to them {\it ex post
facto}, i.e.\ it must be sampled from a conditional distribution.

Coupling into and from the past (see \fref{binary-backoff-ciaftp})
is to coupling into and into the past
(\fref{ciaitp}) as binary-backoff coupling from
the past algorithm (\fref{binary-backoff-cftp}) is to the coupling
into the past (\fref{citp}).  The reference Markov chain
is still run backwards in time, but to test for coalescence, the
random maps of the target Markov chain are composed going forwards in
time.  In this way, as before, the algorithm need only maintain the
images of the random maps of the target chain as the maps are
composed.  Observe that the state of the reference Markov chain at any
given time contains implicit information about the random mappings of
the target Markov chain at all previous times.  This implicit
information can be taken into account when determining the possible
states of the target Markov chain at time 0.  Making use of this
implicit information about previous not-yet-generated random maps is
what distinguishes coupling into and from the past from ordinary CFTP,
and enables it to generate perfectly random samples using
``non-uniformly ergodic'' Markov chains, which cannot be done using
ordinary CFTP.

\begin{figure}[htb]
\begin{code}
\item  $X[0] :=$ ReferenceChainRandomState()
\item  $T := 1$
\item  \Repeat \{
\item\s  SetRandomSeed(seed1[$\log_2(T)$])
\item\s  \For $t := \lfloor T/2 \rfloor+1$ \To $T$
\item\s\s  $X[t]$ := ReverseReferenceChain($X[t-1]$)
\item\s  Set := $\langle$portion of state space compatible with $X[T]\rangle$
\item\s  \For $t := T$ \DownTo $1$
\item\s\s  \If $t$ is a power of $2$
\item\s\s\s  SetRandomSeed(seed2[$\log_2(t)$])
\item\s\s  ApplyTargetChainRandomMapCoupledExPostFacto($X[t]$,$X[t-1]$,Set)
\item\s  $T := 2*T$
\item  \} \Until Singleton(Set)
\item  \Output ElementContainedIn(Set)
\end{code}
\caption{
  Pseudocode for the binary-backoff implementation of CIAFTP.
  The states $\ldots,X_{-t},\ldots,X_{-1},X_0$ of the reference Markov
  chain are generated going back in time and stored, and then read
  going forwards in time; \protect\cite{fill:interruptible} describes
  a way to store many fewer values without excessive recomputation.
  Since the reference and target Markov chains are coupled,
  knowledge of $X[T]$ contains information about previous not-yet-generated random maps of the target chain, which may rule out some values of the state of the target Markov chain.
  The variable Set represents the image of the
  composition $f_{-T}\circ\cdots\circ f_{-t}$ when restricted to this set of possible values, or more generally a
  superset of the image, and is almost never a naive listing of the states.}
\label{fig:binary-backoff-ciaftp}
\end{figure}

{\it Remark:\/} Since (1) the coupling into and from the past
algorithm accesses the random maps of the target chain going forwards
in time, (2) these random maps are coupled {\it ex post facto\/} to
the sample path of the chain, and (3) the sample path of the reference
chain is generated going backwards in time, it is necessary to either
store in memory the entire sample path of the reference chain, or else to
regenerate portions of it as needed.  A similar situation exists in
Fill's algorithm, and \cite{fill:interruptible} describes how to store
portions of the sample path so that not too much memory is used, yet
so that not too much time is spent regenerating the path.

{\it Remark:\/} A few years ago it was asserted that CFTP could not be
used with the so-called non-uniformly ergodic Markov chains.  However,
a variety of algorithms based on coupling into and from the past
(e.g.\ \cite{kendall:area-interaction},
\cite{kendall-moller:exact-spatial}, and \cite{lund-wilson:storage}) do
in fact generate perfectly random samples using non-uniformly ergodic
Markov chains.  One possible interpretation of this asserted
impossibility result is that coupling into and from the past contains
within it another idea.  We are optimistic that there may be
additional clever ideas in the area of perfect simulation.

\section{Read-once coupling into and from the past}
\label{sec:rociaftp}

Since the random maps used in coupling into and from the past are not
independent of one another, but rather are generated by a Markov
process (the reference chain), condition 1 of the main result is not
satisfied.  Therefore our main result does not imply that coupling
into and from the past can be run with a read-once source of
randomness.  In this section we give a protocol for read-once coupling
into and from the past.

For read-once coupling into and from the past, as with read-once CFTP,
it is convenient to work with a composite random map that has an
approximately 1/2 chance of being coalescent.  Of course now the
composite random map takes as input in initial state for the reference
Markov chain, and produces a final state for the reference Markov
chain as well as a random map for the target Markov chain.  The
probability of coalescence is in general a function of the initial
state $X$ for the reference chain, but when this initial state is
random (drawn from the steady state distribution $\pi_{\text{ref}}$ of
the reference Markov chain), the composite random map that we will use
has probability 1/2 of being coalescent.

Imagine that we have a sequence of such composite maps $f_t$ defined
for each integer time $t$.  Let $C_t$ be the indicator random variable
which is $1$ if the $t$th composite map $f_t$ is coalescent, and $0$
otherwise.  In the case of read-once CFTP, the random variables $C_t$
are i.i.d.  All that we can say here is that the $C_t$ form a
stationary process, which is to say that we may shift their indices by
one, and obtain an identically distributed process.

Let $A$ be the set of states $X$ of the reference Markov chain for
which the composite map procedure, when given $X$ as input, generates
a map that is coalescent with probability at least 1/4.  Markov's
inequality tells us that a random $X$ (drawn from $\pi_{\text{ref}}$)
will with probability at least 1/3 be contained in the set $A$.  Since
the reference Markov chain is ergodic, the composite map procedure
projected down to the state space of the reference Markov chain is
also ergodic, so that regardless of its current state, if we wait long
enough (where ``long enough'' may depend upon the current state), with
probability at least 1/6 the state will lie in the set $A$.  The
probability of the next composite map being coalescent is thus at
least 1/24.  So we see that with probability 1, there is some positive
$t$ such that $C_t=1$.  Since the $C_t$'s are stationary process, with
probability 1 there is also some negative $t$ such that $C_t=1$.

Let $S$ be the smallest positive number such that $C_{-S}=1$, and let
$T$ be the smallest non-negative number such that $C_T=1$.  More
generally, let $S_k$ and $T_k$ be the values that $S$ and $T$ would
assume if the indices of the $C_t$ process were decremented by $k$.
Thus the gap between the two coalescent composite maps straddling the
state that time $k$ is given by $S_k+T_k$.  By stationarity, each pair
$(S_k,T_k)$ is distributed in the same manner as $(S,T)$.  Consider
these variables as we increase $k$.  Unless $T_k=0$, we have that
$T_{k+1}=T_k-1$ and $S_{k+1}=S_k+1$.  If on the other hand $T_k=0$,
then $T_{k+1}$ takes on some new non-negative value and $S_{k+1}=1$.
By increasing $k$ many times and averaging, we see that given $S+T$,
the pair $(S,T)$ is uniformly distributed amongst (legal) pairs whose
sum is $S+T$.  In particular, $S$ has the same distribution as $T+1$.

To determine the state of the target Markov chain at time $0$, we can
determine the smallest positive number $S$ such that $f_{-S}$ is
coalescent, and output $f_{-1}\circ f_{-2}\circ \cdots\circ f_{-S}$.
Alternatively, we can randomly generated $S$ from its appropriate
distribution, then using new randomness, generate a sequence of $S$
random composite maps conditioned so that the first is coalescent
while the remaining ones are not, and return the composition of these
maps.  Since $S$ has the same distribution as $T+1$ above, we may
generate random composite maps until one of them is coalescent, and
set $S$ to be the number of generated maps.  To generate the sequence
of $S$ composite maps conditioned so that only the first of them
(rather than the last of them) is coalescent, we can use rejection
sampling.  At no point do we need to reuse previously used randomness.
Pseudocode for this procedure is given in \fref{read-once-ciaftp}.

\begin{figure}[phtb]
\begin{code}
\item  RejectionSample(Length)
\item\s  $X :=$ ReferenceChainRandomState()
\item\s  State := $\langle\text{arbitrary state compatible with $X$}\rangle$
\item\s  RandomComposite($X$,State,CoalescenceFlag)
\item\s  \If \Not CoalescenceFlag
\item\s\s  \Return RejectionSample(Length)
\item\s  \For Count:=2 \To Length 
\item\s\s  RandomComposite($X$,State,CoalescenceFlag)
\item\s\s  \If CoalescenceFlag
\item\s\s\s  \Return RejectionSample(Length)
\item\s  \Return State
\item  
\item  ReadOnceCIAFTP()
\item\s  $X$ := ReferenceChainRandomState()
\item\s  State := $\langle\text{arbitrary state compatible with $X$}\rangle$
\item\s  Count := 0
\item\s  \Repeat
\item\s\s  RandomComposite($X$,State,CoalescenceFlag)
\item\s\s  Count := Count$+1$
\item\s  \Until CoalescenceFlag
\item\s  \Return RejectionSample(Count)
\end{code}
\caption{
  Pseudocode for read-once coupling into and from the past.  Here the
  random composite map is assumed to be defined by an external
  procedure, since it is called from multiple locations.  (Recall that
  the composite map was implicitly defined in the pseudocode for
  read-once coupling from the past, since it would have been called
  from only one location.)  The Set2 argument to RandomComposite is
  optional.  While this procedure does work, we don't recommend it, on
  account of its infinite expected running time.  We leave as a open
  problem the task of finding a more reasonable (finite expected
  running time) read-once version of coupling into and from the past.
  }
\label{fig:read-once-ciaftp}
\end{figure}

Next we consider the expected running time of this read-once coupling
into and from the past procedure.  Let $p_k=\Pr[S=k]$.  When
RejectionSample($k$) it is called, the expected number of trials
before one is accepted is $1/p_k$. Since $T+1$ and $S$ are
equidistributed, with probability $p_k$ ReadOnceCIAFTP() calls
RejectionSample($k$).  Thus the expected number times that
RejectionSample() gets called is $\sum_k p_k \times 1/p_k = \infty$,
unless of course only finitely many $p_k$'s are nonzero.

Despite the expected running time being infinite, since $\sum_k p_k =
1$, the algorithm does terminate with probability 1.

{\it Remark:\/} Upon reading a preliminary explanation of how to do
read-once CFTP, Duncan Murdoch suggested a version that involves the
doubling of starting times in the past used by the binary-backoff CFTP
protocol.  This version fairs poorly when compared with the read-once
CFTP protocol given in \fref{read-once-cftp}.  But when the
binary-backoff variation of read-once CFTP is adapted to coupling into
and from the past, since the protocol in \fref{read-once-ciaftp} has
infinite expected running time, there appears to be no reason to
prefer either variation to the other.  There are many other variations
that work as well, but it is not clear whether or not there is a
variation that has finite expected running time.

\section{Locally stable point processes}
\label{sec:point-process}

The principal application to which coupling into and from the past has
been applied are the locally stable point processes considered by
\cite{kendall-moller:exact-spatial}.  While we don't have a general-purpose read-once
coupling into and from the past protocol, or at least not one that
runs in finite expected time, in this section we see how to apply
read-once CFTP to generate samples from many locally stable point
processes within a reasonable amount of time.

To apply read-once CFTP we need to construct a suitable composite
random map, which is coalescent fairly frequently, and which does not
take any state information as input, such as the value of the
dominating Markov chain, as the auxiliary state information would
introduce dependencies between subsequent random maps.  Since the
natural Markov chain is typically not uniformly ergodic (which was the
reason for introducing the dominating chain in the first place), we
use Murdoch's technique of mixing the natural Markov chain with an
independence sampler.

Following Kendall and M\o ller's notation, we let $f$ denote the
probability density function of the point process relative to a
Poisson process with parameter $\lambda$ (where $\lambda$ is often
constant throughout the region, but can in general vary within the
region).  The process is called ``locally stable'' if for some $K$
adding a point to any given configuration increases its density by no
more than a factor of $K$.
The natural Markov
chain that \cite{kendall-moller:exact-spatial} use is the continuous-time
birth-and-death chain in which points die with rate $1$, and are
proposed with rate $K\lambda$.  A proposed point $x$ is born, i.e.\
added to the configuration $\sigma$, with probability given by
$f(\sigma\cup\{x\}) / [K f(\sigma)]$.  Since detailed balance is satisfied,
this Markov chain has the desired stationary distribution.

Assume for the time being that the density function is always
positive, as it is for say the Strauss process.  (For the Strauss
process, $f(\sigma) = \gamma^{\text{\# pairs of points closer than
$r$}}$ for some interaction radius $r$ and parameter $\gamma<1$
\citep{strauss:strauss} \citep{kelly-ripley:strauss}.)
Later we explain how to deal with point processes such as the
impenetrable spheres model where the density function can be zero.

Our representation for a set of configurations of the point process
will consist of an integer $k$ and two finite sets of points $\Delta$
and $L$.  The set represented by $(k,\Delta,L)$ consists of all
configurations that contain each point in $L$, possibly some of the
points in $\Delta$, and at most $k$ additional points that could be
anywhere.

The first step of the composite map is to generate a Poisson point
configuration with suitably high intensity parameter, such as
$2 K \lambda$, and use it as a proposal for a Metropolis-Hastings update.
Say that the proposed state $\sigma$ has $\#\sigma$ points.  For typical
states of the point process, the proposal $\sigma$ will likely be
rejected, but for all starting states with at least
$B(\sigma) = \#\sigma + \#\sigma \log_2 K + \log_2 1/f(\sigma)$ points,
the proposal will be accepted.  After the
update we know that the state has at most $B(\sigma)$
points, but the points could be anywhere.  Thus the (super)set of possible
configurations is represented by $(B(\sigma),\emptyset,\emptyset)$.

Next we proceed to do updates according to the usual birth-and-death
process.  The $k$ points at unknown locations each die with rate $1$,
so the integer in the representation is decremented at rate $k$.
Points are proposed with rate $K\lambda$, which are then added to
$\Delta$.  The points in $\Delta$ also die with rate $1$.  The set $L$
remains empty.  Eventually $k=0$, implying that each of the original points at
unknown locations has died, and so the Markov chain
state must be a subset of $\Delta$.  Since we now have a finite upper bound
and lower bound on the current state, we can continue to run the
birth-and-death process and start updating the representation of
possible configurations as described by
\cite{kendall-moller:exact-spatial}, i.e.\ using for instance monotone
or anti-monotone coupling.  Eventually $\Delta=\emptyset$, at which
time the the Markov chain is guaranteed to be in state $L$.  Let $T$
be the amount of time that the birth-and-death process was run, i.e.\ 
the time between the Metropolis-Hastings update and the time that
coalescence is achieved.  We throw out the coalesced state and
remember $T$.

Then we do essentially the same thing again, except that we run the
birth-and-death process for an amount of time equal to $T$ rather than
until the time that coalescence is achieved.  We need to evaluate the
random map at the input state as well as detect whether or not it is
coalescent.  Let $\sigma'$ be the Metropolis-Hastings proposal during
this second round, and let $\tau_{\text{old}}$ be the state after
applying Metropolis-Hastings update to the input.  When running the
birth-and-death process during this second round, we maintain an
integer $j$ and finite sets of points $\tau_{\text{old}}$,
$\tau_{\text{new}}$, $\Delta$, and $L$.  All these sets (except
$\tau_{\text{old}}$) are initialized to $\emptyset$.  The current
state resulting from the input state is $\tau = \tau_{\text{old}} \cup
\tau_{\text{new}}$.  When running the birth-and-death process,
any point that is born into $\tau$ is inserted into
$\tau_{\text{new}}$, but of course points in both sets die with rate
$1$.  The integer $k$ is given by $k=j+\#\tau_{\text{old}}$, so we
initialize $j$ to $B(\sigma')-\#\tau_{\text{old}}$ and decrement $j$
at rate $j$.  Thus we can
detect coalescence in the same manner as we did in the first round,
and whether or not coalescence is achieved the second time, we can
determine the final state for any given starting state.  As usual, the
random map preserves the desired probability distribution since both
the Metropolis-Hastings update and the birth-and-death process
preserve the distribution, and because the birth-and-death process is
run for an amount of time that is independent of the birth-and-death
process itself.  If coalescence is in fact officially achieved the
second time, then the coalescence flag is set to be true; with
probability $1/2$ this flag is set to be true.  Thus the composite map
satisfies all the requisite properties to be used with read-once CFTP.

This composite map procedure for locally stable point processes is one
reason that in \sref{unbounded} we recommended doing only the first
update using the independence sampler, and doing subsequent updates
with the natural Markov chain.  While it may be {\it possible\/} to
use the independence sampler more frequently here, doing so would at
the very least unnecessarily complicate the procedure.

\fref{strauss} shows perfectly random samples drawn from
the Strauss point process and the impenetrable spheres model,
which were generated using the approach described here.

\begin{figure}[phtb]
\centerline{\psfig{figure=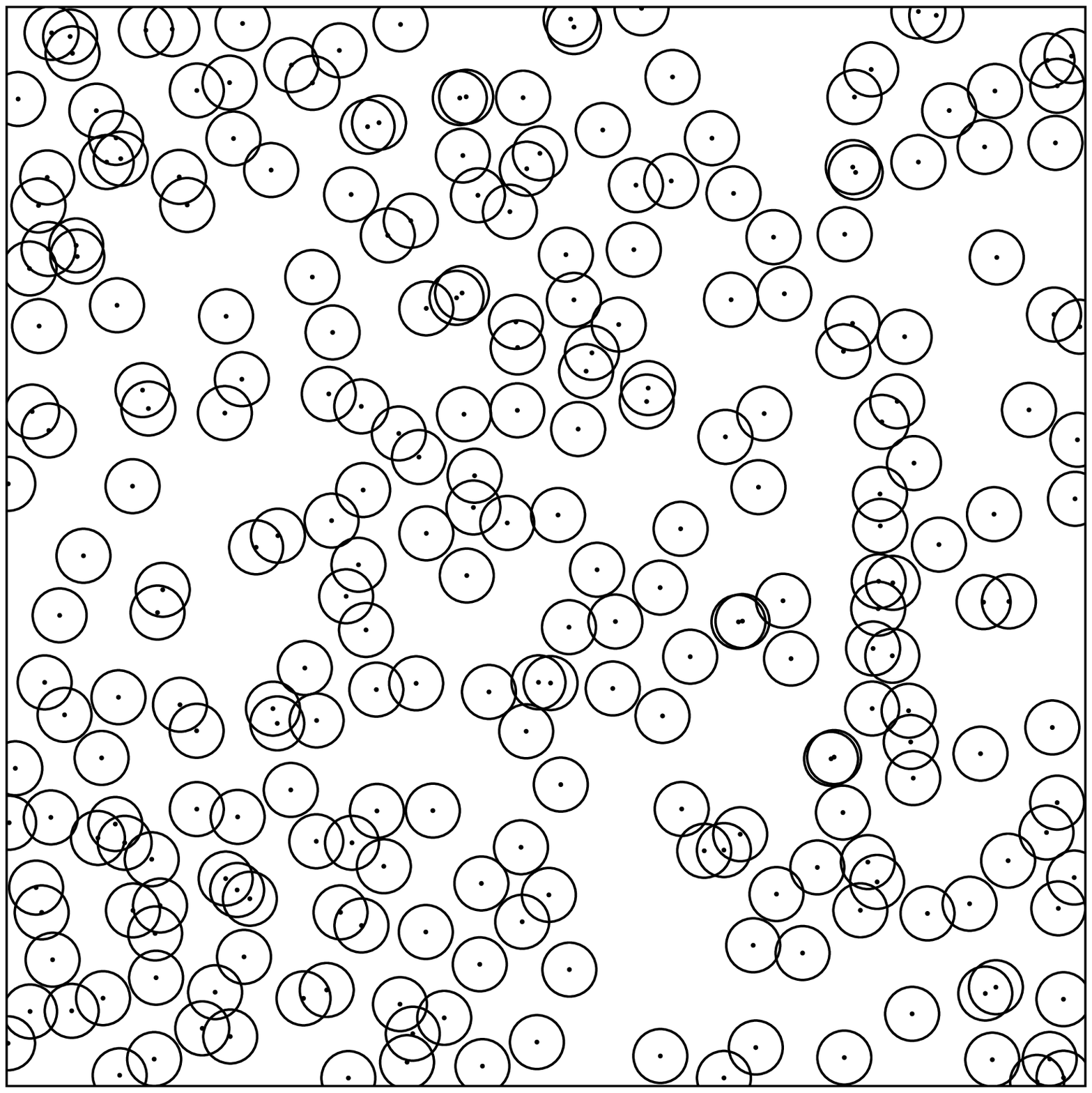,width=3in}\hfil\psfig{figure=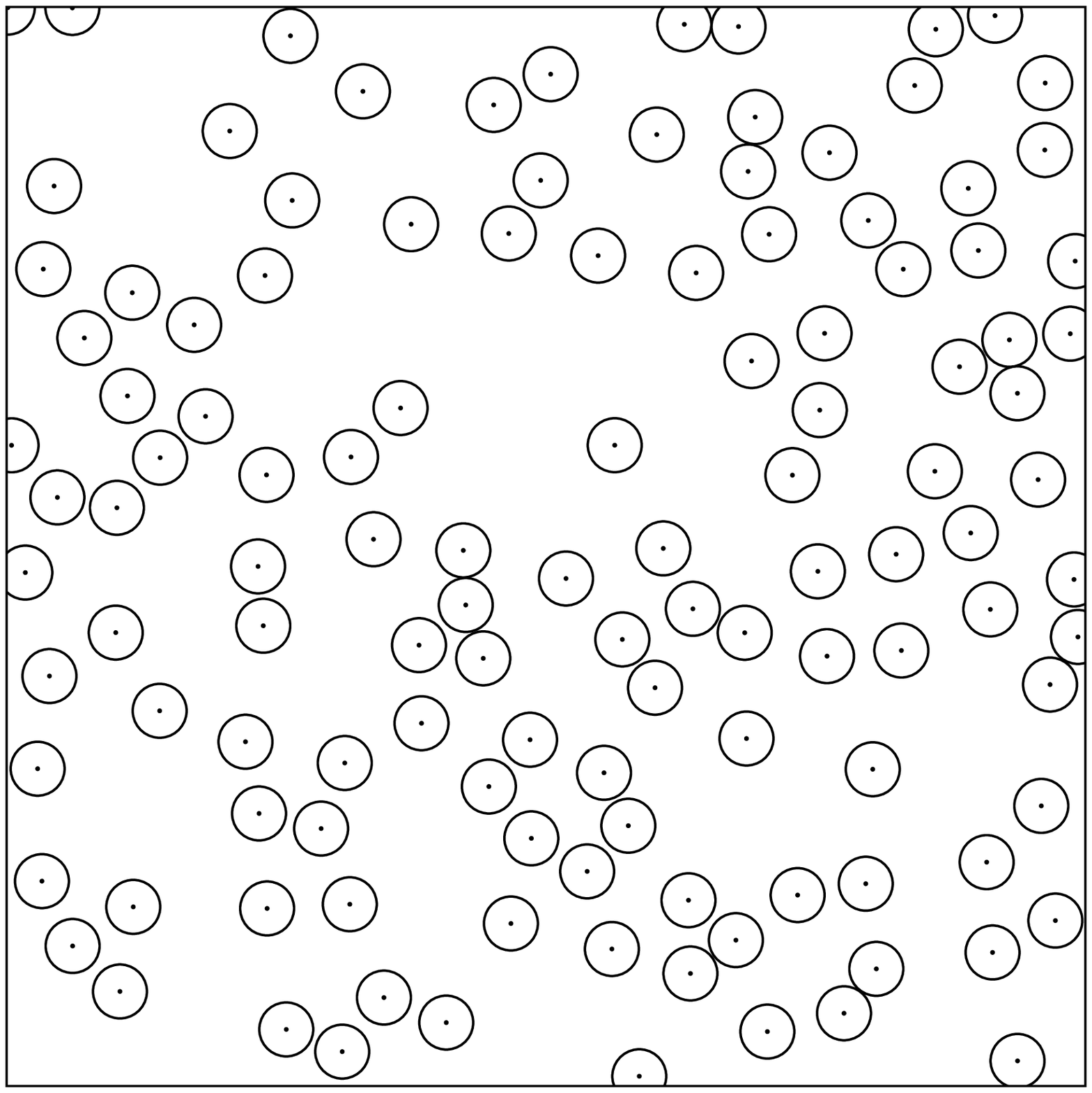,width=3in}}
\caption{
Perfectly random samples of the Strauss point process.  In both panels the points have interaction radius $1$ and occur in a $20\times20$ region with free boundary conditions.  Two points interact, contributing a factor of $\gamma$ to the probability density, if the circles drawn around them overlap.  On the left $\lambda=2$ and $\gamma=1/2$.  On the right $\lambda=1$ and $\gamma=0$ (the impenetrable spheres model).  When $\gamma=0$, it is an immediate consequence of \protect\citep{luby-vigoda:hard-core-glauber} that the coupling is rapid for $\lambda<2/\pi$; in practice it is rapid for $\lambda\leq 1$.
}
\label{fig:strauss}
\end{figure}

Next we consider the distribution of the coalescence time $T$, which
we can break apart as $T=T_1+T_2$, where $T_1$ is the time for all the
points at unknown locations to die, and $T_2$ is the remaining time
for coalescence.  At time $T_1$, the set $\Delta$ is stochastically
dominated by a Poisson point process with intensity $K\lambda$.  Since
\cite{kendall-moller:exact-spatial} coalesce all the states that are
subsets of such a Poisson point process, the remaining time $T_2$
until coalescence will be at most as large as the coalescence time for
Kendall and M\o ller.  The time $T_1$ will be about $\log B(\sigma)+O(1)$.
One might object that the time it takes the procedure to decrement $k$
to $0$ is not proportional to $T_1$, but is instead proportional its
initial value $B(\sigma)$.
But since $T_1$ distributed in the same manner as $-\log(1-U^{1/B(\sigma)})$,
where $U$ is uniformly distributed between 0 and 1, we can optimize
the procedure to avoid needless decrementing.
Recall that $B(\sigma) = \#\sigma + \#\sigma \log_2 K + \log_2 1/f(\sigma)$
where $\sigma$ is a Poisson point process with intensity $2 K \lambda$.
Thus we can expect $T_1$ to be fairly small unless $f(\sigma)$ is
extraordinaraly close to $0$.  As $f(\sigma)\rightarrow 0$, the time
$T_1$ grows like $\log \log 1/f(\sigma)$.

This very weak dependence of the time $T_1$ upon $f(\sigma)$ allows us
to run the algorithm even for point processes such as the impenetrable
spheres model where the density function $f(\sigma)$ can be zero.
What we do is ``soften'' the density, and work with a modified
density.  For the case of the impenetrable spheres model, rather than
multiply the density function by $0$ for each pair of points that are
too close together, we can instead multiply the density by say
$10^{-20}$ for each such pair.  We generate point configurations from
this modified density, and then do rejection sampling to obtain a
sample from the desired density.  Since $10^{-20}$ is extremely close
to $0$, we almost never reject the sample, and since the run time
dependence upon $f(\sigma)$ is so weak, the run time is not adversely
affected.

In practice the time $T_2$ is either on par with or else dwarfs the
time $T_1$.  Since the run time of the algorithm used by
\cite{kendall-moller:exact-spatial} is only marginally larger than
$T_2$, we expect the running time of their algorithm and ours to be
similar, with smaller memory requirements for the algorithm described
here.  Our algorithm for the Strauss process also appears to be
competitve with another algorithm given by
\cite{moller-nicholls:perfect-tempering}.

\section{Summary and open problems}
\label{sec:open}

We have given the modification of the coupling from the past protocol
which only requires a read once source of randomness.  This read-once
CFTP protocol is on par with the usual CFTP protocol in terms of
memory and time, and for some applications will be up to
logarithmically faster.  Read-once CFTP is closely related to the
PASTA property from operations research.  We have also given a
read-once version of the coupling into and from the past protocol, but
it is unsatisfactory since the expected running time is infinite.  We
leave as open problems the existence of a better read-once version of
CIAFTP, and the existence of further connections between PASTA-type
thereoms and perfect sampling algorithms.

\section*{Source code}

The source code for the program used to make the Strauss process
samples in \fref{strauss} is available at
\texttt{http://dbwilson.com/strauss/}.

\section*{Acknowledgements}

The author thanks Wilfrid Kendall, Duncan Murdoch, and Jim Propp
for their comments on earlier versions of this manuscript.

\bibliography{exact,tree,rocftp}
\bibliographystyle{plainnat}

\end{document}